\newcommand{\Z}{{\mathbb Z}}
\newcommand{\Q}{{\mathbb Q}} 
 
\newcommand{\R}{{\mathbb R}}
\newcommand{\RZ}{{\mathbb R_0}}
\newcommand{\RO}{{\mathbb R_1}}
\newcommand{\NZ}{{\mathbb N_0}}
\newcommand{\NO}{{\mathbb N_1}}
\newcommand{\QZ}{{\mathbb Q_0}}
\newcommand{\QO}{{\mathbb Q_1}}
\newcommand{\QT}{{\Q\left[(2)\right]}}

\newcommand{\BU}{{\widetilde{U}}}
\newcommand{\BOmega}{{\widetilde{\Omega}}}
\newcommand{\xto}{{\;\xrightarrow{\,\boldsymbol +}\;}}

\documentclass[oneside,final,12pt,leqno]{amsart}
\usepackage{amsmath}        
\usepackage{amsthm}         
\usepackage{latexsym}       
\usepackage{amssymb}

\title[The Real $3x+1$ Problem]{The Real $\mathbf{3x+1}$ Problem}

\author{Pavlos B. Konstadinidis}
\address{Departamento de Matem\'atica, \hfill\break\indent  
Instituto de Matem\'atica e Estat\'\i stica, \hfill\break\indent  
Universidade de S\~ao Paulo \hfill\break\indent 
C.P.G. Sala B-23 Rua do Mat\~ao, 1010 \hfill\break\indent
S\~ao Paulo, S.P., Brazil 05508-090 \hfill\break\indent}
\email{pavlos@ime.usp.br}
\thanks{The author was supported by CNPq-Brazil}
\subjclass[2000]{11B37, 26A18}

\theoremstyle{plain}\newtheorem{teo}{Theorem}[section]
\theoremstyle{plain}\newtheorem{lem}[teo]{Lemma}
\theoremstyle{plain}\newtheorem{prop}[teo]{Proposition}
\theoremstyle{plain}
\theoremstyle{definition}
\theoremstyle{remark}\newtheorem{rem}[teo]{Remark}

\begin{document}
\maketitle

\begin{section}{Introduction}\label{sec:intro}
\indent First of all, let's fix some notations. \\
\indent As usual, $\R$, $\Q$ and $\Z$ will denote the sets of all real,
rational and integer numbers, respectively. \\
\indent Put $\RO = \{x\in \R : x\ge 1\}$, $\,\,\QO = \Q\cap\RO\,$ and
$\,\,\NO = \Z\cap\RO = \{1,2,3,\ldots\}$. \\
\indent Put $\RZ = \{x\in \R : x\ge 0\}$, $\,\,\QZ = \Q\cap\RZ\,$ and 
$\,\,\NZ = \Z\cap\RZ = \{0,1,2,\ldots\}$. \\
\indent For $x\in\R$, $\lfloor x\rfloor$ will denote the {\em floor} or
{\em integer part} of $x$, that is to say, $\lfloor x\rfloor
=\,$max$\,\{k\in\Z: k\le x\}$. \\
\\
\indent The well-known $3n+1$ function (see, e.g.,~\cite{LAG3} and~\cite{WIR})
is the function $T:\NO\to\NO$ given by
\begin{equation}
T(n)=\left\{ \begin{array}{lll}
               T_0(n)={\,\,\,\,\,\mbox{\Large{$\frac{n}{2}$}}} & \mbox{ if $n$ is even,} \\
      \\
               T_1(n)={\mbox{\Large{$\frac{3n+1}{2}$}}}    & \mbox{ if $n$ is odd.}
               \end{array} \right.
\end{equation}
\indent In this work, we introduce another extension of $T$, namely the function $U:\RO\to\RO$ defined by
\begin{equation}
\,\,\,\,\,U(x)=\left\{ \begin{array}{lll}
               U_0(x)={\,\,\,\,\,\mbox{\Large{$\frac{x}{2}$}}} & \mbox{ if
               $\lfloor x\rfloor$ is even,} \\
      \\
               U_1(x)={\mbox{\Large{$\frac{3x+1}{2}$}}}    & \mbox{ if $\lfloor
               x\rfloor$ is odd.}
               \end{array} \right.
\end{equation}
\indent Note that $U\vert_{\textrm{\small{$\NO$}}}$ (the restriction of $U$ to
$\NO$) is indeed $T$. We shall call $U$ the {\em real} $3x+1$ function 
(in contrast to the {\em integer} $3n+1$ function $T$). In
Section~\ref{sec:conj}, we'll propose a conjecture about the iterates of $U$
that generalizes the famous $3n+1$ conjecture. We'll then prove our main 
result about the iterates of $U$ (Theorem~\ref{thm:cycles}), which is directly 
related to both of these conjectures. We'll also introduce the {\em flipped} 
$3x+1$ function $\BU$ and prove an analogous result for its iterates. In
Section~\ref{sec:obs}, we'll show a couple of simple propositions about the
iterates of $U$ and $\BU$, introduce other related functions and propose some
questions and conjectures about their iterates. \\
\indent We hope that the results, 
conjectures and questions stated here will be not only relevant to the 
$3n+1$ conjecture itself, but also of interest in their own right. All of the 
results presented here were independently discovered by the author,
but some of them already appear in the literature. In those cases, we refer
the reader to their proofs. However, for our and the reader's benefit, we do 
recall some well-known definitions (in a format slightly better suited
to our purposes). 
\end{section}
\begin{section}{The Conjecture and the Main Results}\label{sec:conj}
\indent Given a (nonempty) set $X$ and a function $f:X\to X$, the
{\em iterates of} $f$ will be denoted by $f^i$ ($i\in \NZ$). They are defined 
by $f^0=$ id$_{X}$ (the identity function on X) and by $f^i=f\circ f^{i-1}$, 
for $i>0$. For any given $x\in X$, the $f$-{\em trajectory of} $x$ or  
{\em starting at} $x$ is the sequence $\mathcal
T_f(x)=\left(f^i(x)\right)_{i=0}^{\infty}$. An 
$f$-{\em periodic trajectory} or, simply, an $f$-{\em cycle} is the 
$f$-trajectory of some $z\in X$ such that $f^n(z)=z$, for some $n\in \NO$ (in
this case, the $f$-cycles starting at $f^k(z)$, $k\in\NZ$, will sometimes be 
considered as being one and the same $f$-cycle). By an $f$-{\em cycle of 
length} $l\in\NO$ we mean any sequence in the set 
$\left\{\left(x,f(x),\ldots,f^l(x)\right):x\in X,\,f^l(x)=x\right\}$.\\
\indent Now, let $\QT$ denote the set of all rational numbers having an odd
denominator when written in lowest terms (see \cite{HAL & HUN}). A number
$a/b\in\QT$ (with an odd $b$) is {\em even} ({\em odd}$\,$) if its numerator
$a$ is even (odd). The {\em rational Collatz sequence generated} by $r_0\in
\QT$ is the $g$-trajectory of $r_0$, where $g:\QT\to \QT$ is given by 
$g(r)=g_0(r)=r/2$, if $r$ is even, and
$g(r)=g_1(r)=(3r+1)/2$, if $r$ is odd. A {\em rational Collatz cycle} ({\em of
length $l\,$}) is simply a $g$-cycle (of length $l$). 
Given $l\in\NO$ and $n\in \NZ$, let $S_{l,n}$ be 
the set of all $0$-$1$ sequences of length $l$ containing exactly $n$ $1$'s, and put $S_l=\bigcup_{n=0}^lS_{l,n}$ and $S=\bigcup_{l=1}^{\infty}S_l$. If   
$s\in S$, we'll denote the number of $1$'s in $s$ by $n(s)$ and the
length of $s$ by $l(s)$. Given $s=(s_1,s_2,\ldots,s_l)\in S$, 
define $\phi_s:\R\to \R$ by $\phi_s=g_{s_l}\circ \cdots \circ g_{s_2}\circ g_{s_1}$. A sequence
$(x_0,x_1,\ldots,x_l)$ of numbers $x_i\in \R$ is called a {\em pseudo-cycle of
  length} $l$ if there exists $s=(s_1,s_2,\ldots,s_l)\in S$ such that
$x_l=x_0$ and $x_i=g_{s_i}(x_{i-1})$, for $i=1,2,\ldots,l$ (note that $x_l=\phi_s(x_0)$). Finally, define $\varphi:S\to \NZ$ by
$\varphi(s)=\sum_{j=1}^{l(s)}s_j2^{j-1}3^{s_{j+1}+s_{j+2}+\cdots+s_{l(s)}}$. \\
\indent Let $x_0\in\RO$ be given. If 
$\left\{\lim_{k\to\infty}U^{2k}(x_0),\,\lim_{k\to\infty}U^{2k+1}(x_0)\right\}=\{1,2\}$, 
then we'll say that its $U$-trajectory $\mathcal T_U(x_0)$ {\em tends to}
$\{1,2\}$ and this will be denoted by $\mathcal T_U(x_0)\to\{1,2\}$. Our 
{\em real $3x+1$ conjecture} is \\
\\
\indent {\bf RU}: For all $x\in\RO$, $\mathcal T_U(x)\to \{1,2\}$. \\
\\
\indent Note that, for all $n\in\NO$, $\mathcal T_T(n)=\mathcal T_U(n)$. The 
famous (integer) $3n+1$ conjecture may then be stated as \\
\\
\indent {\bf NU}: For all $n\in\NO$, $\mathcal T_U(n)\to \{1,2\}$. \\
\\
\indent One could also state both of these conjectures in terms of the $U$-{\em parity sequence} associated with $x\in\RO$, which is simply
the infinite $0$-$1$ sequence 
$\mathcal P_U(x)=\left(\left\lfloor U^i(x)\right\rfloor\bmod 2\right)_{i=0}^{\infty}$. 
Note that this
sequence encodes which branch of $U$ ($U_0$ or $U_1$) is used in each step of
$\mathcal T_U(x)$. Now, an infinite $0$-$1$ sequence $(p_i)_{i=0}^{\infty}$ 
will be called {\em eventually periodic with period} $(0,1)$ if there exists
$j\in\NZ$ such that $(p_{i},\,p_{i+1})=(0,1)$, for all $i=j+2m$, $m\in\NZ$. 
It's a simple matter (see Proposition~\ref{prop:period}) to show that, 
for each $x\in\RO$, $\mathcal P_U(x)$ is eventually periodic with period
$(0,1)$, if, and only if, $\mathcal T_U(x)\to\{1,2\}$. In other words, 
the conjectures {\bf RU} and {\bf NU} above can be stated in the following 
alternative, equivalent forms. \\
\\
\indent {\bf RU'}: For all $x\in\RO$, $\mathcal P_U(x)$ is eventually periodic
with period $(0,1)$. \\
\indent {\bf NU'}: For all $n\in\NO$, $\mathcal P_U(n)$ is eventually periodic
with period $(0,1)$. \\
\\
\indent Now, we observe that our ${\bf RU}$ conjecture clearly implies both
of the following two conjectures. \\
\\
\indent {\bf OU}: The only $U$-cycle is the trivial $T$-cycle $(1,2,1,2,1,\ldots)$. \\
\indent {\bf BU}: Every $U$-trajectory is bounded. \\
\\
\indent Of course, all $T$-cycles are $U$-cycles, and one would naturally 
expect to find (many) more $U$-cycles than $T$-cycles. However, 
our main result, which is directly related to the conjectures 
{\bf RU} and {\bf OU} above, tells us that in fact quite the 
opposite happens. 
\begin{teo}\label{thm:cycles}
The only $U$-cycles are the $T$-cycles. 
\end{teo}
\begin{proof} 
\indent Let's first state next two lemmas that will be used in this and
subsequent proofs. The reader may find their proofs in~\cite{HAL & HUN}
and~\cite{LAG2} (the basic idea of most of Lemma~\ref{lem:rat} below
is due originally to B\"ohm and Sontacchi~\cite{BOHM}). 

\begin{lem}\label{lem:rat}
{\em (B\"ohm and Sontacchi, Lagarias, Halbeisen and Hungerb\"uhler).} A sequence
$(x_0,x_1,\ldots,x_l)$ is a rational Collatz cycle of length $l$ if, and only
if, it is a pseudo-cycle of length $l$. Moreover, if a rational Collatz cycle
is not the cycle $(0,0,\ldots)$, then it's elements are either all strictly 
positive or all strictly negative. $\Box$ 
\end{lem}

\begin{lem}\label{lem:coll}
{\em (Lagarias).} For any $s\in S$ and any $x\in\R$, we have
 that 
\begin{equation} \phi_s(x)=\frac{3^{n(s)}x+\varphi(s)}{2^{l(s)}} \mbox{.}
\end{equation}
\indent Therefore, given $s\in S$, 
\begin{equation} x_0(s)=\frac{\varphi(s)}{2^{l(s)}-3^{n(s)}}\in\QT 
\end{equation}
is the {\bf unique} number that generates the rational Collatz cycle of
length $l(s)$ that is also the pseudo-cycle of length $l(s)$ determined by
$s$. $\Box$ 
\end{lem}

\indent To begin with, we note thal all $U$-cycles start at numbers in 
$\QO$, since, for each $k\in\NO$, we have that every solution of $x=U^k(x)$ is 
rational. Let's suppose then that there exist $x_0\in\QO\setminus\NO$ and 
$l\in\NO$ such that there's a $U$-cycle of length
$l$ starting at $x_0$, namely
$\Omega(x_0)=\left(x_0,U(x_0),\ldots,U^l(x_0)=x_0\right)$. If we derive a 
contradiction from this hypothesis, then we'll be done. Note that it's
immediate (by inspection) that the only $U$-cycle of length less than $4$ is
the $T$-cycle $(1,2,1)$. Hence, without loss of generality, we may assume that
$l\ge 4$, which avoids our having to treat some trivial cases separately in 
what follows. Now, since $U_\iota\equiv g_\iota$ ($\iota=0,1$), 
$\Omega(x_0)$ is a pseudo-cycle of length $l$. Thus, by Lemma~\ref{lem:rat}, 
$\Omega(x_0)$ is a rational Collatz cycle of length $l$ as well. Therefore, by
using Lemma~\ref{lem:coll} and the fact that 
$U\vert_{\textrm{\small{$\NO$}}}=T$, one obtains both that
all $U^i(x_0)\in\QT\cap\QO\setminus\NO$ and that  
\begin{equation} x_0=x_0(s)=\frac{\varphi(s)}{2^{l(s)}-3^{n(s)}}\mbox{,} 
\end{equation}
where $s=(s_1,s_2,\ldots,s_l)\in S$ is the $0$-$1$ sequence associated with
(the pseudo-cycle) $\Omega(x_0)$, i.e., $s$ consists of the first $l=l(s)$
terms in $\mathcal P_U(x_0)$. For convenience, put $n=n(s)$ and $d=2^l-3^n$. Now, given
any $a/b$ in $\QT$ (with an odd $b$), it's clear that every term in the
rational Collatz sequence generated by $a/b$ may be written with denominator
$b$. As $d$ happens to be odd, one may, for $i=0,1,\ldots,l$, write that   
\begin{equation}
  x_i=U^i(x_0)=\frac{c_i}{d}=\frac{q_id+r_i}{d}=q_i+\frac{r_i}{d}\mbox{,} 
\end{equation}
where $q_i$ is the quotient and $r_i$ the remainder in the Euclidean division
of $c_i$ by $d$. Note that all $c_i,\,q_i$ and $r_i$ lie in
$\NO$ and that $d\geq 5$ (for all $x_i=c_i/d$ are in $\QO\setminus\NO$,
$\varphi(s)>0$ and 3 doesn't divide $d$). In particular, no $r_i$ is $0$, and 
so all $r_i$ satisfy $0<r_i<d$. Moreover, because $d=2^l-3^n>0$, one has
\begin{equation} n<l\log_32.\label{eq:ineq0}
\end{equation}
\indent Now, since
$\Omega(x_0)=(x_0,x_1,\ldots,x_l=x_0)$ is both a $U$-cycle and a rational
Collatz cycle (of length $l$), we have, for $i=0,1,\ldots,l$, that
$q_i=\left\lfloor x_i\right\rfloor$ is even (odd) if, and only if, $c_i$ is
even (odd). Thus, all $r_i=c_i-dq_i$ are even. Write $r_i=2^{e_i}o_i$, where
$e_i\ge 1$ and $o_i$ is odd, think of $r_0,r_1,\ldots,r_l=r_0$ as being
arranged (in this order) in a circular manner and observe that, 
for $i=0,1,\ldots,l$,  
\begin{equation} r_i= \left\{ \begin{array}{lllll}
               \,\,\,\,\,\,\mbox{{\Large{$\frac{1}{2}$}}$\,r_{i-1}$} & \mbox{ if $q_{i-1}$ is even,} \\
     \\
               \,\,\,\,\,\,\mbox{{\Large{$\frac{3}{2}$}}$\,r_{i-1}$} & \mbox{ if $q_{i-1}$
is odd and 
               $r_{i-1}<$\,{\Large{$\frac{2}{3}$}}$\,d$,} \\
      \\
               \mbox{{\Large{$\frac{3}{2}$}}$\,r_{i-1}-d$} & \mbox{ if $q_{i-1}$ is odd and 
               $r_{i-1}>$\,{\Large{$\frac{2}{3}$}}$\,d$.} 
               \end{array} \right. \label{eq:rs}
\end{equation}
Note that, since $3$ doesn't divide $d$, it's never the case that
$r_{i-1}=2d/3$ in (\ref{eq:rs}). As usual, indices are to be considered modulo 
$l$ whenever it's the case to do so. Now, if $r_i$ is such that
$r_i=3r_{i-1}/2-d$, then we'll say that this $r_i$ is {\em
  new}. Note that, if $r_{j+1}$ is not new, we have that either
$r_{j+1}=2^{e_j-1}o_j$ or $r_{j+1}=2^{e_j-1}(3o_j)$. This clearly means that
at least one of $r_0,r_1,\ldots,r_{l-1}$ is new. By renaming the
$x_i$'s if necessary, we can assume that $r_0$ ($=r_l$) is new. Now, let 
$0\le p<q\le l$ be such that $r_p$ and $r_q$ are {\em consecutive new}, 
that is, both $r_p$ and $r_q$ are new and, for all $p<k<q$, $\,r_k$ is not new 
(if $r_0$ is the only new one, then put $p=0$ and $q=l$). Because there're no
new $r_k$'s strictly between $r_p$ and $r_q$, one has 
\begin{equation} r_{q-1}=2\left(3^{n(p,\,q-1)}o_p\right)
  \,\,\,\mbox{and}\,\,\,\, e_p=q-p \mbox{,} \label{eq:eq0}
\end{equation}
where, for any $0\le i\le j\le l$, $n(i,j)$ is the number of times $U_1$ is
used from $x_i$ to $x_j$, i.e., $n(i,j)$ is the number of $1$'s in
$\{s_{i+1},s_{i+2},\ldots,s_j\}$. Since $r_q$ is new, we have both that
$n(p,q)=n(p,q-1)+1$ and that $d<3r_{q-1}/2$. From this and (\ref{eq:eq0}), it
follows that 
\begin{equation}
d<3^{n(p,\,q-1)+1}o_p=\frac{3^{n(p,\,q-1)+1}}{2^{e_p}}\,r_p=\frac{3^{n(p,\,q)}}{2^{q-p}}\,r_p
\mbox{.} \label{eq:ineq1}
\end{equation} 
\indent Now, because $r_p$ is new, $r_p=3r_{p-1}/2-d$, and so, since
$0<r_{p-1}<d$, we obtain $r_p<d/2<2d/3$. From this and (\ref{eq:ineq1}), one gets 
\begin{equation}
\indent d<\frac{3^{n(p,\,q)-1}}{2^{q-p-1}}\,d\Longrightarrow
  3^{n(p,\,q)-1}>2^{q-p-1}\Longrightarrow
  n(p,q)>\log_32^{q-p-1}+1 \mbox{.} 
\end{equation}
\indent Therefore, $n(p,q)>\log_32^{q-p-1}+\log_32=\log_32^{q-p}$, and so one
has that 
\begin{equation} n(p,q)>(q-p)\log_32 \mbox{.} \label{eq:ineq2}
\end{equation}
\indent Now, let $0=i_0<i_1<\cdots <i_m=l$, $m\ge 1$, be such that 
$r_{i_0},r_{i_1},\ldots,r_{i_m}$ are all the new $r_i$'s in
$\{r_0,r_1,\ldots,r_l\}$. We have that $n=\sum_{k=1}^mn(i_{k-1},\,i_k)$,
$l=\sum_{k=1}^m(i_k-i_{k-1})$ and that $r_{i_{k-1}}$ and $r_{i_k}$ are
consecutive new for all $k=1,2,\ldots,m$. Consequently, inequality
(\ref{eq:ineq2}) gives us $n>l\log_32$, but this last inequality 
{\bf contradicts} inequality (\ref{eq:ineq0}). 
\end{proof}
\indent We note that some authors have already investigated a variety 
of interesting {\em smooth} 
extensions of $T$ to the real (and even complex) numbers 
(see, e.g.,~\cite{CHAM1},~\cite{CHAM2},~\cite{DUM},~\cite{LAG1} 
and~\cite{LETH}). Unlike the {\em conjectured} case of $U$, 
however, the dynamics of these extensions outside the integers are 
always extraneous to the $3x+1$ conjecture (i.e., there exist periodic 
and divergent trajectories). \\ 
\indent Now, the previous theorem illustrated the relative ease one has in obtaining
some results if he is allowed the freedom to work in $\RO$ (instead 
of his having to concentrate solely on $\NO$). For another example along these 
lines, consider the {\em flipped} $3x+1$ function $\BU:\RZ\to\RZ$ defined by 
\begin{equation}
\,\,\,\,\,\BU(x)=\left\{ \begin{array}{lll}
               \BU_0(x)=U_1(x) & \mbox{ if
               $\lfloor x\rfloor$ is even,} \\
      \\
               \BU_1(x)=U_0(x) & \mbox{ if $\lfloor
               x\rfloor$ is odd.}
               \end{array} \right.
\end{equation}
\indent Clearly, $\BU\vert_{\textrm{\small{$\NZ$}}}$ 
is not a function from $\NZ$ to $\NZ$. Naturally, one would like to know  
what happens to the $\BU$-trajectories. In particular, one would try to obtain 
all $\BU$-cycles. This is in fact done in our next theorem, which is a 
bonus result we've gotten from the method we've used to prove Theorem~\ref{thm:cycles}. 

\begin{teo}\label{thm:nobar}
There are no $\BU$-cycles.
\end{teo}
\begin{proof}
The proof is almost entirely analogous to the proof of Theorem~\ref{thm:cycles}
above, and so we'll be brief and point out only the required modifications. 
Clearly, no $\BU$-cycles start at numbers in $\NZ$. Let's suppose then that 
there exist $x_0\in\QZ\setminus\NZ$ and $l\in\NO$ such that there's a
$\BU$-cycle of length $l$ starting at $x_0$, namely
$\BOmega(x_0)=(x_0,\BU(x_0),\ldots,\BU^l(x_0)=x_0)$. If this assumption leads
us to a contradiction, then we'll be done. By inspection, there're no 
$\BU$-cycles of length less than $4$, and so we may, without loss of
generality, assume that $l\geq 4$ (again, this assumption is made so as to
avoid trivialities in what follows). Now, with similar notations and the same 
arguments from the proof of Theorem~\ref{thm:cycles}, one obtains, for
$i=0,1,\ldots,l$, that  
\begin{equation}
  x_i=\BU^i(x_0)=\frac{c_i}{d}=\frac{q_id+r_i}{d}=q_i+\frac{r_i}{d}=(q_i+1)-\frac{d-r_i}{d}\mbox{,} 
\end{equation}
where $q_i$ is the quotient and $r_i$ the remainder in the Euclidean division
of $c_i$ by $d$. Since no $x_i$'s belong to $\NZ$, we have that all 
$r_i$ satisfy $0<d-r_i<d$. Moreover, because $d=2^l-3^n>0$, we have, as before,
that  
\begin{equation} n<l\log_32.\label{eq:ineqbar}
\end{equation}
\indent Since $\BOmega(x_0)=(x_0,x_1,\ldots,x_l=x_0)$ is both a
$\BU$-cycle and a rational Collatz cycle (of length $l$), it follows, for
$i=0,1,\ldots,l$, that $q_i=\left\lfloor x_i\right\rfloor$ is even (odd) if,
and only if, $c_i$ is odd (even). Thus, all $r_i=c_i-dq_i$ are odd, i.e., all
$d-r_i$ are even. Now, think of $d-r_0,d-r_1,\ldots,d-r_l=d-r_0$ as being 
arranged (in this order) in a circular fashion and note that, for 
$i=0,1,\ldots,l$, 
\begin{equation} d-r_i= \left\{ \begin{array}{lllll}
               \,\,\,\,\,\,\mbox{{\Large{$\frac{1}{2}$}}$\,(d-r_{i-1}$)} & \mbox{ if $q_{i-1}$ is odd,} \\
     \\
               \,\,\,\,\,\,\mbox{{\Large{$\frac{3}{2}$}}$\,(d-r_{i-1}$)} &
     \mbox{ if $q_{i-1}$
is even and 
               $d-r_{i-1}<$\,{\Large{$\frac{2}{3}$}}$\,d$,} \\
      \\
               \mbox{{\Large{$\frac{3}{2}$}}$\,(d-r_{i-1})-d$} & \mbox{ if
$q_{i-1}$ is even and 
               $d-r_{i-1}>$\,{\Large{$\frac{2}{3}$}}$\,d$.} 
               \end{array} \right. \label{eq:rsbar} 
\end{equation}
\indent Now, arguing exactly in the same way as we've done in the proof of 
Theorem~\ref{thm:cycles}, we conclude that $n>l\log_3 2$, which contradicts 
(\ref{eq:ineqbar}). 
\end{proof}
\indent Note that yet another equivalent way of phrasing the  
conjecture {\bf RU} is to say that, for every $x\in\RO$, there exists
$k\in\NZ$ such that $U^k(x)\in [1,3)$. Our corresponding conjecture for the 
iterates of $\BU$ is \\
\\
\indent {\bf R}{$\mathbf\BU$}: For every $x\in\RZ$ there exists $k\in\NZ$ such
that $\BU^k(x)\in [0,2)$. \\
\\
\indent Of course, Theorem~\ref{thm:nobar} is directly related to
the conjecture {\bf R}{$\mathbf\BU$} above. Let's conclude this section
by observing that our {\bf R}{$\mathbf\BU$} conjecture clearly implies the 
following conjecture. \\
\\
\indent {\bf B}{$\mathbf\BU$}: Every $\BU$-trajectory is bounded.
\end{section}

\begin{section}{Other Results, Conjectures and Questions}\label{sec:obs}
One way to find out if studying what happens to the iterates of $U$ can shed 
some new light on the $3n+1$ conjecture or not would be to try and answer our 
first question. \\
\\
{\bf Q1}: Does the $3n+1$ conjecture imply our real $3x+1$ conjecture {\bf RU}? \\
\\
\indent On one hand, if the answer to this question is {\em yes}, then this
would show that looking at the iterates of $U$ amounts to essentially the same
thing as looking at those of $T$ (as far as the $3n+1$ conjecture is
concerned). On the other hand, we note that, if the $3n+1$ conjecture is
true, then the answer to the question {\bf Q1} above could very well be 
{\em no}. To see how this might be so, suppose that, instead of $T$, one
considered the original Collatz function, i.e., the function 
$f:\NO\to\NO$ given by 
\begin{equation}
f(n)=\left\{ \begin{array}{lll}
               f_0(n)={\,\,\,\,\,\mbox{\Large{$\frac{n}{2}$}}} & \mbox{ if $n$ is even,} \\
      \\
               f_1(n)={\mbox{$3n+1$}}    & \mbox{ if $n$ is odd.}
               \end{array} \right.
\end{equation} 

Its extension to $\RO$ (in our sense) is the function $F:\RO\to\RO$ given by
\begin{equation}
F(x)=\left\{ \begin{array}{lll}
               F_0(x)={\,\,\,\,\,\mbox{\Large{$\frac{x}{2}$}}} & \mbox{ if
               $\lfloor x\rfloor$ is even,} \\
      \\
               F_1(x)={\mbox{$3x+1$}}    & \mbox{ if $\lfloor x\rfloor$ is odd.}
               \end{array} \right.
\end{equation} 

The statement for the $F$-trajectories which corresponds to the  
conjecture {\bf RU} would be the claim that, for all $x\in\RO$, 
$\mathcal T_F(x)\to\{1,2\}$. However, this is readily seen to be
{\em false}, since one has, for example, that all $F$-trajectories starting at
$2m+3/2$, $m\in\NZ$, diverge (monotonically) to $+\infty$. Now, the 
$3n+1$ conjecture for the iterates of $T$ is equivalent to the (same) one for
the iterates of $f$. 
Thus, if the $3n+1$ conjecture turns out to be true, then the 
question for the $F$-trajectories that is the counterpart to 
question {\bf Q1} will have a negative answer. Moreover, if our real $3x+1$ 
conjecture {\bf RU} is true, then the $U$-trajectories and the
$F$-trajectories will be seen to have quite different behaviors in $\RO$ 
(as opposed to what happens in $\NO$). In our view, comparisons between
the $U$-trajectories and the $F$-trajectories may play an important 
r\^ole in some future $3x+1$-type investigations. Let our next question 
emphasize this point. \\
\\
{\bf Q2}: Are the $F$-trajectories starting at $2m+3/2$, $m=0,1,2,\ldots$, 
the only $F$-trajectories that do not tend to $\{1,2\}$? \\
\\
\indent Of course, analogous questions on similar notions regarding the 
iterates of $\BU$ could be posed as well. We'll now show a simple result about
the iterates of $U$. Its proof will suggest a new approach one might consider 
in trying to prove the conjecture {\bf OU} (see Remark~\ref{rem:ou}). 
A corresponding result for the iterates of $\BU$ will be then obtained as a corollary. 
Before we can state these results, a couple of definitions are needed. \\
\indent Given $x_0\in\RO$, we'll say that
$\mathcal P_U(x_0)=(p_i)_{i=0}^{\infty}$ is {\em eventually periodic with
  period} $s=(s_0,s_1,\ldots,s_{l(s)-1})\in S$ if there exists $j\in\NZ$ such
that $(p_i,p_{i+1},\ldots,p_{i+l(s)-1})=(s_0,s_1,\ldots,s_{l(s)-1})$ for all
$i=j+ml(s)$, $m\in\NZ$. Moreover, if $a\in\NO$ is such that there's
a $U$-cycle of length $l$ starting at $a$, then we'll say that $\mathcal
T_U(x_0)$ {\em tends to} $\{U^t(a)\}$ {\em from above} (in symbols,
$\mathcal T_U(x_0)\xto\{U^t(a)\}$) if there is $j_0\in\NZ$ such that, 
for all $j\in\{0,1,\ldots,l-1\}$, $U^{kl}\left(U^{j+j_0}(x_0)\right)\to U^j(a)^{\boldsymbol +}$, as $k\to+\infty$. 
\begin{prop}\label{prop:period}
If $a\in\NO$ is such that there's a $U$-cycle of length $l$ starting at
$a$, then, for all $x\in\RO$, we have that $\mathcal P_U(x)$ is eventually
periodic with period $(a\bmod 2,U(a)\bmod 2,\ldots,U^{l-1}(a)\bmod 2)$, if, and
only if, $\mathcal T_U(x)$ tends to $\{U^t(a)\}$ from above. 
\end{prop}
\begin{proof}
Suppose at first that $x\in\RO$ is such that $\mathcal
T_U(x)\xto\{U^t(a)\}$. Now, from the fact that there's a 
$U$-cycle of length $l$ starting at $a$, it clearly follows that there's some 
$0<\theta\in\R$ such that, for all $y\in [a,a+\theta)$ and all $m\in\NZ$, 
\begin{equation}
\left(\lfloor U^{ml}(y)\rfloor,\lfloor
U^{ml+1}(y)\rfloor,\ldots,\lfloor
U^{ml+l-1}(y)\rfloor\right)=\left(a,U(a),\ldots,U^{l-1}(a)\right).  
\end{equation}
For instance, any $0<\theta<(2/3)^l$ will do. Since $\mathcal
T_U(x)\xto\{U^t(a)\}$, there's some $k_0\in\NZ$ such that
$U^{k_0}(x)\in [a,a+\theta)$. Hence, $\mathcal P_U(x)$ is eventually
periodic with period $(a\bmod 2,U(a)\bmod 2,\ldots,U^{l-1}(a)\bmod 2)$. 
For the other direction, suppose now that $x\in\RO$ is such that $\mathcal
P_U(x)$ is eventually periodic with period 
$s=(a\bmod 2,U(a)\bmod 2,\ldots,U^{l-1}(a)\bmod 2)\in S$. 
By using Lemma~\ref{lem:coll}, one sees that there's some $j_0\in\NZ$ such
that 
\begin{multline*}
U^{l(s)}\left(U^{j_0}(x)\right)=\frac{3^{n(s)}U^{j_0}(x)+\varphi(s)}{2^{l(s)}}=\frac{3^{n(s)}\left(a+U^{j_0}(x)-a\right)+\varphi(s)}{2^{l(s)}}=\\[6 pt] 
=\frac{3^{n(s)}a+\varphi(s)}{2^{l(s)}}+\frac{3^{n(s)}\left(U^{j_0}(x)-a\right)}{2^{l(s)}}=a+\frac{3^{n(s)}}{2^{l(s)}}\left(U^{j_0}(x)-a\right).
\end{multline*}
Analogously, we have, for $j=0,1,\ldots,l-1$, that 
\begin{equation*}
U^{l(s)}\left(U^{j+j_0}(x)\right)=U^j(a)+\frac{3^{n(s)}}{2^{l(s)}}\left(U^{j+j_0}(x)-U^j(a)\right).
\end{equation*}
Therefore, for all $m\in\NZ$ and all $j\in\{0,1,\ldots,l-1\}$, 
\begin{equation}
U^{ml(s)}\left(U^{j+j_0}(x)\right)=U^j(a)+\left(\frac{3^{n(s)}}{2^{l(s)}}\right)^{\!\!m}\!\!\left(U^{j+j_0}(x)-U^j(a)\right).
\end{equation}
Since $3^{n(s)}<2^{l(s)}$, it's not hard to conclude now that 
$\mathcal T_U(x)\xto\{U^t(a)\}$.
\end{proof}
\indent Now, with the appropriate analogous definitions for the iterates of
$\BU$, the same argument presented in the proof of
Proposition~\ref{prop:period} above gives us the following result as well. 
\begin{prop}
If $a\in\NO$ is such that there's a $U$-cycle of length $l$ starting at
$a$, then, for all $x\in\RZ$, we have that 
$\mathcal P_\BU(x)$ is eventually periodic with period 
$\left(1-a\bmod 2,1-\left(U(a)\bmod 2\right),\ldots,1-\left(U^{l-1}(a)\bmod
    2\right)\right)$ if, and only if, $\mathcal T_\BU(x)$ tends to $\{U^t(a)\}$ from below. $\Box$
\end{prop}

\indent Note that if the $3n+1$ conjecture is true and $x_0\in\RO$ is such 
that $\mathcal T_U(x_0)\to a_0^{\boldsymbol +}$ for some $a_0\in\NO$, 
then $\mathcal T_U(x_0)\to\{1,2\}$. This indicates one way in which one may 
try and give a positive answer to question {\bf Q1}. \\
\indent Now, consider $\mathcal I_U(\NO)=\{x\in\RO:\mathbf{\exists}\,k\in\NZ\textrm{ with } U^k(x)\in\NO\}$ and $\mathcal N_U(\NO)=\RO\setminus\,\mathcal I_U(\NO)$. Of course, our {\bf RU} conjecture implies the following conjecture. \\
\\
{\bf $\mathcal N$U}: For all $x\in\mathcal N_U(\NO)$, $\mathcal T_U(x)\to\{1,2\}$. \\
\\
\indent We may pose now our next question, which can also be thought of as 
being one of the possible (non-trivial) ways of turning question {\bf Q1}
around. \\
\\
{\bf Q3}: Does the {\bf $\mathcal N$U} conjecture above imply the $3n+1$ conjecture? 

\begin{rem}\label{rem:ou}
Let's just note here an interesting corollary of the proof of
Proposition~\ref{prop:period}: if one proves that, for all $n\in\NO$ and 
all $0<\rho\in\R$, there exists some $z\in (n,n+\rho)\cap\mathcal 
N_U(\NO)$ such that $\mathcal T_U(z)\to\{1,2\}$, then it will follow that 
the {\bf OU} conjecture (which, in light of Theorem~\ref{thm:cycles}, is in
fact the ``there are no non-trivial $T$-cycles'' conjecture) is true. 
\end{rem}

\indent To try and answer the question {\bf Q3} above might be an  
even better way of seeing whether there're some real advantages in 
shifting one's attention from $T$ to $U$. Let's end this line of 
inquiries now by our registering the following very broad (but also
potentially very productive) question. \\
\\ 
{\bf Q4}: What kind of results for the iterates of $U$ does one get by 
attempting to translate known results for the iterates of $T$? \\
\\ 
\indent In conclusion, let's just remark that the apparent general project 
would be for one to study the dynamical system in $\R$ generated by
the iterates of the (discontinuous) piecewise linear functions of the 
following ``simple'' kind. \\
\indent Let $\alpha,\beta,\gamma,\delta,\tau\in\R$ be fixed,
with $\tau\in [0,2)$, and consider the function
$\Phi=\Phi(\alpha,\beta,\gamma,\delta,\tau):\R\to\R$ defined by  
\begin{equation}\label{eq:lambda}
\\\Phi(x)=\left\{ \begin{array}{lll}
               \Phi_0(x)={\mbox{$\alpha x+\beta$}} & \mbox{
               if $\lfloor x+\tau\rfloor$ is even,} \\
      \\
               \Phi_1(x)={\mbox{$\gamma x+\delta$}}    & \mbox{ if
               $\lfloor x+\tau\rfloor$ is odd.}
               \end{array} \right.
\end{equation}
\indent Naturally, the crux of the matter here is to find out how the
parameters $\alpha,\beta,\gamma,\delta$ and $\tau$ affect the behavior of
the $\Phi=\Phi(\alpha,\beta,\gamma,\delta,\tau)$-trajectories. This
brings us to our final (albeit seemingly intractable as of yet!) question. \\
\\
{\bf Q5}: How do the general properties of the dynamical
system in $\R$ generated by the iterates of the function 
$\Phi=\Phi(\alpha,\beta,\gamma,\delta,\tau)$ defined as in
(\ref{eq:lambda}) depend on the values of the real parameters 
$\alpha,\beta,\gamma,\delta$ and $\tau$? \\
\\
\indent E.g., $U=\Phi(1/2,0,3/2,1/2,0)\vert_{\textrm{\small{$\RO$}}}$ and 
$\BU=\Phi(1/2,0,3/2,1/2,1)\vert_{\textrm{\small{$\RZ$}}}$. Note also that 
the functions $\Phi(1/2,0,3/2,1/2,\tau_0)\vert_{\textrm{\small{$\RO$}}}$, 
with $0\leq\tau_0<1$, are all extensions of $T$. Finally, we bring into 
attention $V=\Phi(1/2,0,3/2,0,0)\vert_{\textrm{\small{$\RO$}}}$, i.e., 
the function $V:\RO\to\RO$ given by
\begin{equation}
V(x)=\left\{ \begin{array}{lll}
               V_0(x)=\,\mbox{{\Large{$\frac{1}{2}$}}$\,x$} & \mbox{ if
               $\lfloor x\rfloor$ is even,} \\
      \\
               V_1(x)=\,\mbox{{\Large{$\frac{3}{2}$}}$\,x$} & \mbox{ if $\lfloor x\rfloor$ is odd.}
               \end{array} \right.
\end{equation} 
\indent Of course, there're no $V$-cycles. It might be worthwhile for one 
to try and find out the status of our following final two conjectures, as well 
as their possible connections to the $3n+1$ and {\bf RU} conjectures, if any. 
\\
\\
{\bf RV}: For every $x\in\RO$ there exists $k\in\NZ$ such that $V^k(x)\in
[1,3)$. \\
{\bf BV}: Every $V$-trajectory is bounded. \\ 
\\ 
\indent The author would like to thank the referee for suggestions that 
have lead to an improvement in the presentation of this paper. 
\end{section}

\end{document}